Т.Г. Эргашев

# ЧЕТВЕРТЫЙ ПОТЕНЦИАЛ ДВОЙНОГО СЛОЯ ДЛЯ ОБОБЩЕННОГО ДВУОСЕСИММЕТРИЧЕСКОГО УРАВНЕНИЯ ГЕЛЬМГОЛЬЦА

Потенциал двойного слоя играет важную роль при решении краевых задач для эллиптических уравнений. При его исследовании существенно используются свойства фундаментальных решений данного уравнения. В настоящее время все фундаментальные решения обобщенного двуосесимметрического уравнения Гельмгольца известны, но, несмотря на это, только для первого из них построена теория потенциала. В данной работе исследуется потенциал двойного слоя, соответствующий четвертому фундаментальному решению. При использовании свойства гипергеометрической функции Аппеля от двух переменных доказываются предельные теоремы и выводятся интегральные уравнения, содержащие в ядре плотности потенциала двойного слоя.

**Ключевые слова:** *обобщенное двуосесимметрическое уравнение Гельмгольца; формула Грина; фундаментальное решение; четвертый потенциал двойного слоя; гипергеометрические функции Аппеля от двух переменных; интегральные уравнения с плотностью потенциала двойного слоя в ядре.*

### 1. Введение

Многочисленные приложения теории потенциала можно найти в механике жидкости, эластодинамике, электромагнетизме и акустике. С помощью теории потенциала краевые задачи удаётся свести к решению интегральных уравнений.

Потенциал двойного слоя играет важную роль при решении краевых задач для эллиптических уравнений. При этом решение ищется в виде потенциала двойного слоя с неизвестной плотностью, для определения которой применяется теория интегральных уравнений Фредгольма второго рода [1–3]. В свою очередь, такой потенциал выписывается через фундаментальное решение данного эллиптического уравнения.

Фундаментальные решения следующего обобщенного двуосесимметрического уравнения Гельмгольца:

$$H_{\alpha,\beta}^{\lambda}(u) \equiv u_{xx} + u_{yy} + \frac{2\alpha}{x}u_x + \frac{2\beta}{y}u_y - \lambda^2 u = 0,$$

здесь $\alpha$, $\beta$ и $\lambda$ – постоянные, причем $0 < 2\alpha, 2\beta < 1$, приведены в [4]. Оказывается, когда $\lambda = 0$, все четыре фундаментальные решения $q_i(x, y; x_0, y_0)$ ($i = 1, 2, 3, 4$) уравнения

$$H_{\alpha,\beta}^{0}(u) \equiv u_{xx} + u_{yy} + \frac{2\alpha}{x}u_x + \frac{2\beta}{y}u_y = 0 \qquad (1.1)$$

можно выразить с помощью гипергеометрической функции Аппеля от двух переменных второго рода $F_2(a, b_1, b_2; c_1, c_2; x, y)$, определенной по формуле [5–7]

$$F_2(a;b_1,b_2;c_1,c_2;x,y) = \sum_{m,n=0}^{\infty} \frac{(a)_{m+n}(b_1)_m(b_2)_n}{(c_1)_m(c_2)_n m!n!} x^m y^n,$$

где $(a)_n$ – символ Похгаммера: $(a)_0 = 1$, $(a)_n = a \cdot (a+1) \cdot (a+2) \cdot ... \cdot (a+n-1)$, $n = 1, 2, 3, ...$

Теория потенциала для простейшего вырождающегося эллиптического уравнения (т.е. при $\alpha = 0$ и $\lambda = 0$) изложена автором [8, 9]. Следуя его теории, в работе [10] для первого фундаментального решения $q_1(x,y;x_0,y_0)$ уравнения (1.1) построена теория потенциала двойного слоя в области

$$\Omega \subset R_+^2 = \{(x,y): x > 0, y > 0\}.$$

В данной работе исследуется потенциал двойного слоя, соответствующий четвертому фундаментальному решению уравнения (1.1):

$$q_4(x,y;x_0,y_0) = k_4 (r^2)^{\alpha+\beta-2} x^{1-2\alpha} y^{1-2\beta} x_0^{1-2\alpha} y_0^{1-2\beta} \times$$
$$\times F_2(2-\alpha-\beta; 1-\alpha, 1-\beta; 2-2\alpha, 2-2\beta; \xi, \eta), \qquad (1.2)$$

где

$$k_4 = \frac{2^{4-2\alpha-2\beta}}{4\pi} \frac{\Gamma(1-\alpha)\Gamma(1-\beta)\Gamma(2-\alpha-\beta)}{\Gamma(2-2\alpha)\Gamma(2-2\beta)},$$

$$\xi = \frac{r^2 - r_1^2}{r^2}, \quad \eta = \frac{r^2 - r_2^2}{r^2}, \quad \left.\begin{array}{c}r^2\\r_1^2\\r_2^2\end{array}\right\} = \left(\begin{array}{c}-\\x+x_0\\-\end{array}\right)^2 + \left(\begin{array}{c}-\\y-y_0\\+\end{array}\right)^2. \qquad (1.3)$$

Нетрудно проверить, что функция $q_4(x,y;x_0,y_0)$ по переменным $(x_0,y_0)$ является решением уравнения (1.1) и обладает следующими свойствами:

$$q_4(x,y;x_0,y_0)\big|_{x=0} = 0, \quad q_4(x,y;x_0,y_0)\big|_{y=0} = 0. \qquad (1.4)$$

Используя свойства гипергеометрической функции Аппеля от двух переменных, доказываем предельные теоремы и выводим интегральные уравнения, содержащие в ядре плотность потенциала двойного слоя.

## 2. Формула Грина

Рассмотрим тождество

$$x^{2\alpha} y^{2\beta} \left[ u H_{\alpha,\beta}^0(v) - v H_{\alpha,\beta}^0(u) \right] = \frac{\partial}{\partial x}\left[ x^{2\alpha} y^{2\beta} (v_x u - v u_x) \right] + \frac{\partial}{\partial y}\left[ x^{2\alpha} y^{2\beta} (v_y u - v u_y) \right].$$

Интегрируя обе части последнего тождества по области $\Omega$, расположенной в первой четверти ($x > 0$, $y > 0$), и пользуясь формулой Остроградского, получим

$$\iint_\Omega x^{2\alpha} y^{2\beta} \left[ u H_{\alpha,\beta}^0(v) - v H_{\alpha,\beta}^0(u) \right] dx dy =$$
$$= \int_S x^{2\alpha} y^{2\beta} u (v_x dy - v_y dx) - x^{2\alpha} y^{2\beta} v (u_x dy - u_y dx), \qquad (2.1)$$

где $S = \partial\Omega$ – контур области $\Omega$.

Формула Грина (2.1) выводится при следующих предположениях: функции $u(x,y)$, $v(x,y)$ и их частные производные первого порядка непрерывны в замкнутой области $\overline{\Omega}$, частные производные второго порядка непрерывны внутри $\Omega$ и интегралы по $\Omega$, содержащие $H^0_{\alpha,\beta}(u)$ и $H^0_{\alpha,\beta}(v)$, имеют смысл. Если $H^0_{\alpha,\beta}(u)$ и $H^0_{\alpha,\beta}(v)$ не обладают непрерывностью вплоть до $S$, то это – несобственные интегралы, которые получаются как пределы по любой последовательности областей $\Omega_n$, которые содержатся внутри $\Omega$, когда эти области $\Omega_n$ стремятся к $\Omega$, так что всякая точка, находящаяся внутри $\Omega$, попадает внутрь областей $\Omega_n$, начиная с некоторого номера $n$.

Если $u$ и $v$ суть решения уравнения (1.1), то из формулы (2.1) имеем

$$\int_S x^{2\alpha} y^{2\beta} \left( u \frac{\partial v}{\partial n} - v \frac{\partial u}{\partial n} \right) ds = 0. \qquad (2.2)$$

Здесь

$$\frac{\partial}{\partial n} = \frac{dy}{ds} \frac{\partial}{\partial x} - \frac{dx}{ds} \frac{\partial}{\partial y}, \quad \frac{dy}{ds} = \cos(n,x), \quad \frac{dx}{ds} = -\cos(n,y), \qquad (2.3)$$

$n$ – внешняя нормаль к кривой $S$.

Полагая в формуле (2.1) $v = 1$ и заменяя $u$ на $u^2$, получим

$$\iint_\Omega x^{2\alpha} y^{2\beta} \left[ u_x^2 + u_y^2 \right] dx dy = \int_S x^{2\alpha} y^{2\beta} u \frac{\partial u}{\partial n} ds,$$

где $u(x,y)$ – решение уравнения (1.1).

Наконец, из формулы (2.2), полагая $v = 1$, будем иметь

$$\int_S x^{2\alpha} y^{2\beta} \frac{\partial u}{\partial n} ds = 0, \qquad (2.4)$$

т.е. интеграл от нормальной производной решения уравнения (1.1) с весом $x^{2\alpha} y^{2\beta}$ по контуру области равен нулю.

## 3. Потенциал двойного слоя $w^{(4)}(x_0, y_0)$

Пусть $\Omega$ – область, ограниченная отрезками $(0,a)$ и $(0,b)$ осей $x$ и $y$ соответственно и кривой $\Gamma$ с концами в точках $A(a,0)$ и $B(0,b)$, лежащей в первой четверти $x > 0,\ y > 0$.

Параметрическое уравнение кривой $\Gamma$ пусть будет $x = x(s),\ y = y(s)$, где $s$ – длина дуги, отсчитываемая от точки $B$. Относительно кривой $\Gamma$ будем предполагать, что:

1) функции $x = x(s)$ и $y = y(s)$ имеют непрерывные производные $x'(s)$ и $y'(s)$ на отрезке $[0,l]$, не обращающиеся одновременно в нуль; вторые произ-водные $x''(s)$ и $y''(s)$ удовлетворяют условию Гельдера на $[0,l]$, где $l$ – длина кривой $\Gamma$;

2) в окрестностях точек $A(a,0)$ и $B(0,b)$ на кривой $\Gamma$ выполняются условия

$$\left|\frac{dx}{ds}\right| \le Cy^{1+\varepsilon}(s), \quad \left|\frac{dy}{ds}\right| \le Cx^{1+\varepsilon}(s), \ 0 < \varepsilon < 1, \tag{3.1}$$

где $C$ – постоянная. Координаты переменной точки на кривой $\Gamma$ будем обозначать через $(x,y)$.

Рассмотрим интеграл

$$w^{(4)}(x_0,y_0) = \int_0^l x^{2\alpha} y^{2\beta} \mu_4(s) \frac{\partial q_4(x,y;x_0,y_0)}{\partial n} ds, \tag{3.2}$$

где $q_4(x,y;x_0,y_0)$ – фундаментальное решение уравнения (1.1), определенное по формуле (1.2), а $\mu_4(s)$ – непрерывная функция в промежутке $[0,l]$.

Интеграл (3.2) будем называть *четвертым потенциалом двойного слоя с плотностью* $\mu_4(s)$. Очевидно, что $w^{(4)}(x_0,y_0)$ есть регулярное решение уравнения (1.1) в любой области, лежащей в первой четверти, не имеющей общих точек ни с кривой $\Gamma$, ни с осью $x$ и ни с осью $y$. Как и в случае логарифмического потенциала, можно показать существование потенциала двойного слоя (3.2) в точках кривой $\Gamma$ для ограниченной плотности $\mu_4(s)$. Потенциал двойного слоя (3.2) при $\mu_4(s) \equiv 1$ обозначим через $w_1^{(4)}(x_0,y_0)$.

**Лемма 1.** *Справедливы следующие формулы:*

$$w_1^{(4)}(x_0,y_0) = \begin{cases} k(x_0,y_0) - 1, & (x_0,y_0) \in \Omega, \\ k(x_0,y_0) - \dfrac{1}{2}, & (x_0,y_0) \in \Gamma, \\ k(x_0,y_0), & (x_0,y_0) \notin \overline{\Omega}, \end{cases} \tag{3.3}$$

*где* $\overline{\Omega} := \Omega \cup \Gamma$;

$$\begin{aligned}
k(x_0,y_0) = (1-2\beta)k_4 x_0^{1-2\alpha} y_0^{1-2\beta} \int_0^a x\left((x-x_0)^2 + y_0^2\right)^{\alpha+\beta-2} \times \\
\times F\left(2-\alpha-\beta, 1-\alpha; 2-2\alpha; -\frac{4xx_0}{(x-x_0)^2 + y_0^2}\right)dx + \\
+ (1-2\alpha)k_4 x_0^{1-2\alpha} y_0^{1-2\beta} \int_0^b y\left(x_0^2 + (y-y_0)^2\right)^{\alpha+\beta-2} \times \\
\times F\left(2-\alpha-\beta, 1-\beta; 2-2\beta; -\frac{4yy_0}{x_0^2 + (y-y_0)^2}\right)dy.
\end{aligned} \tag{3.4}$$

*Здесь* $F(a,b;c;z) = \sum_{k=0}^{\infty} \dfrac{(a)_k (b)_k}{(c)_k k!} z^k$ – *известная гипергеометрическая функция Гаусса.*

***Доказательство. Случай 1.*** Пусть точка $(x_0, y_0)$ находится внутри $\Omega$. Вырежем из области $\Omega$ круг малого радиуса $\rho$ с центром в точке $(x_0, y_0)$ и обозначим через $\Omega^\rho$ оставшуюся часть области $\Omega$, а через $C_\rho$ окружность вырезанного круга. В области $\Omega^\rho$ функция $q_4(x,y;x_0,y_0)$ – регулярное решение уравнения (1.1). Используя следующую формулу для производной гипергеометрической функции Аппеля [4]:

$$\frac{\partial^{m+n} F_2(a;b_1,b_2;c_1,c_2;x,y)}{\partial x^m \partial y^n} =$$
$$= \frac{(a)_{m+n}(b_1)_m(b_2)_n}{(c_1)_m(c_2)_n} F_2(a+m+n;b_1+m,b_2+n;c_1+m,c_2+n;x,y), \quad (3.5)$$

имеем

$$\frac{\partial q_4(x,y;x_0,y_0)}{\partial x} = (1-2\alpha)k_4 (r^2)^{\alpha+\beta-2} x^{-2\alpha} y^{1-2\beta} x_0^{1-2\alpha} y_0^{1-2\beta} \times$$
$$\times F_2(2-\alpha-\beta;1-\alpha,1-\beta;2-2\alpha,2-2\beta;\xi,\eta) -$$
$$-2(2-\alpha-\beta)k_4 (r^2)^{\alpha+\beta-3} x^{1-2\alpha} y^{1-2\beta} x_0^{1-2\alpha} y_0^{1-2\beta} P(x,y;x_0,y_0), \quad (3.6)$$

где

$$P(x,y;x_0,y_0) = (x-x_0)F_2(2-\alpha-\beta;1-\alpha,1-\beta;2-2\alpha,2-2\beta;\xi,\eta) +$$
$$+ x_0 F_2(3-\alpha-\beta;2-\alpha,1-\beta;3-2\alpha,2-2\beta;\xi,\eta) +$$
$$+ (x-x_0)\left[\frac{(1-\alpha)}{2-2\alpha}\xi F_2(3-\alpha-\beta;2-\alpha,1-\beta;3-2\alpha,2-2\beta;\xi,\eta) + \right.$$
$$\left. + \frac{(1-\beta)}{2-2\beta}\eta F_2(3-\alpha-\beta;1-\alpha,2-\beta;2-2\alpha,3-2\beta;\xi,\eta)\right]. \quad (3.7)$$

Далее применяя известное соотношение [4]:

$$\frac{b_1}{c_1}x F_2(a+1;b_1+1,b_2;c_1+1,c_2;x,y) + \frac{b_2}{c_2}y F_2(a+1;b_1,b_2+1;c_1,c_2+1;x,y) =$$
$$= F_2(a+1;b_1,b_2;c_1,c_2;x,y) - F_2(a;b_1,b_2;c_1,c_2;x,y)$$

к квадратной скобке в (3.7), получаем

$$\frac{\partial q_4(x,y;x_0,y_0)}{\partial x} = (1-2\alpha)k_4 (r^2)^{\alpha+\beta-2} x^{-2\alpha} y^{1-2\beta} x_0^{1-2\alpha} y_0^{1-2\beta} \times$$
$$\times F_2(2-\alpha-\beta;1-\alpha,1-\beta;2-2\alpha,2-2\beta;\xi,\eta) -$$

$$-2(2-\alpha-\beta)k_4 (r^2)^{\alpha+\beta-3} x^{1-2\alpha} y^{1-2\beta} x_0^{2-2\alpha} y_0^{1-2\beta} \times$$
$$\times F_2(3-\alpha-\beta;2-\alpha,1-\beta;3-2\alpha,2-2\beta;\xi,\eta) -$$

$$-2(2-\alpha-\beta)(x-x_0)k_4 (r^2)^{\alpha+\beta-3} x^{1-2\alpha} y^{1-2\beta} x_0^{1-2\alpha} y_0^{1-2\beta} \times$$
$$\times F_2(3-\alpha-\beta;1-\alpha,1-\beta;2-2\alpha,2-2\beta;\xi,\eta). \quad (3.8)$$

Аналогично находим

$$\frac{\partial q_4(x,y;x_0,y_0)}{\partial y} = (1-2\beta)k_4\left(r^2\right)^{\alpha+\beta-2} x^{1-2\alpha} y^{-2\beta} x_0^{1-2\alpha} y_0^{1-2\beta} \times$$
$$\times F_2(2-\alpha-\beta;1-\alpha,1-\beta;2-2\alpha,2-2\beta;\xi,\eta) -$$
$$-2(2-\alpha-\beta)k_4\left(r^2\right)^{\alpha+\beta-3} x^{1-2\alpha} y^{1-2\beta} x_0^{1-2\alpha} y_0^{2-2\beta} \times$$
$$\times F_2(3-\alpha-\beta;1-\alpha,2-\beta;2-2\alpha,3-2\beta;\xi,\eta) -$$
$$-2(2-\alpha-\beta)(y-y_0)k_4\left(r^2\right)^{\alpha+\beta-3} x^{1-2\alpha} y^{1-2\beta} x_0^{1-2\alpha} y_0^{1-2\beta} \times$$
$$\times F_2(3-\alpha-\beta;1-\alpha,1-\beta;2-2\alpha,2-2\beta;\xi,\eta). \quad (3.9)$$

Пользуясь (3.8) и (3.9), в силу (1.2) и (2.3) найдем

$$\frac{\partial q_4(x,y;x_0,y_0)}{\partial n} = -(2-\alpha-\beta)k_4\left(r^2\right)^{\alpha+\beta-2} x^{1-2\alpha} y^{1-2\beta} x_0^{1-2\alpha} y_0^{1-2\beta} \times$$
$$\times F_2(3-\alpha-\beta;1-\alpha,1-\beta;2-2\alpha,2-2\beta;\xi,\eta)\frac{\partial}{\partial n}\left[\ln r^2\right] -$$
$$-2(2-\alpha-\beta)k_4\left(r^2\right)^{\alpha+\beta-3} x^{1-2\alpha} y^{1-2\beta} x_0^{2-2\alpha} y_0^{1-2\beta} \times$$
$$\times F_2(3-\alpha-\beta;2-\alpha,1-\beta;3-2\alpha,2-2\beta;\xi,\eta)\frac{dy(s)}{ds} +$$
$$+2(2-\alpha-\beta)k_4\left(r^2\right)^{\alpha+\beta-3} x^{1-2\alpha} y^{1-2\beta} x_0^{1-2\alpha} y_0^{2-2\beta} \times$$
$$\times F_2(3-\alpha-\beta;1-\alpha,2-\beta;2-2\alpha,3-2\beta;\xi,\eta)\frac{dx(s)}{ds} +$$
$$+(1-2\alpha)k_4\left(r^2\right)^{\alpha+\beta-2} x^{-2\alpha} y^{1-2\beta} x_0^{1-2\alpha} y_0^{1-2\beta} \times$$
$$\times F_2(2-\alpha-\beta;1-\alpha,1-\beta;2-2\alpha,2-2\beta;\xi,\eta)\frac{dy(s)}{ds} -$$
$$-(1-2\beta)k_4\left(r^2\right)^{\alpha+\beta-2} x^{1-2\alpha} y^{-2\beta} x_0^{1-2\alpha} y_0^{1-2\beta} \times$$
$$\times F_2(2-\alpha-\beta;1-\alpha,1-\beta;2-2\alpha,2-2\beta;\xi,\eta)\frac{dx(s)}{ds}. \quad (3.10)$$

Тогда, в силу (2.4) и (1.4) получим

$$w_1^{(4)}(x_0,y_0) = \lim_{\rho\to 0}\int_{C_\rho} x^{2\alpha} y^{2\beta} \frac{\partial q_4(x,y;x_0,y_0)}{\partial n} ds +$$
$$+\int_0^b y^{2\beta}\left[x^{2\alpha}\frac{\partial q_4(x,y;x_0,y_0)}{\partial n}\right]_{x=0} ds + \int_0^a x^{2\alpha}\left[y^{2\beta}\frac{\partial q_4(x,y;x_0,y_0)}{\partial n}\right]_{y=0} ds. \quad (3.11)$$

Подставив (3.10) в (3.11), найдем

$$w_1^{(4)}(x_0,y_0) = k_4 x_0^{1-2\alpha} y_0^{1-2\beta} \times$$
$$\times \lim_{\rho\to 0}\left[(2-\alpha-\beta)\{-J_1 - 2x_0 J_2 + 2y_0 J_3\} + (1-2\alpha)J_4 - (1-2\beta)J_5\right] + J_6 + J_7, \quad (3.12)$$

где

$$J_1(x_0, y_0) = \int\limits_{C_\rho} xy\left(r^2\right)^{\alpha+\beta-2} F_2(3-\alpha-\beta; 1-\alpha, 1-\beta; 2-2\alpha, 2-2\beta; \xi, \eta) \frac{\partial}{\partial n}\left[\ln r^2\right] ds;$$

$$J_2(x_0, y_0) = \int\limits_{C_\rho} xy\left(r^2\right)^{\alpha+\beta-3} F_2(3-\alpha-\beta; 2-\alpha, 1-\beta; 3-2\alpha, 2-2\beta; \xi, \eta) \frac{dy(s)}{ds} ds;$$

$$J_3(x_0, y_0) = \int\limits_{C_\rho} xy\left(r^2\right)^{\alpha+\beta-3} F_2(3-\alpha-\beta; 1-\alpha, 2-\beta; 2-2\alpha, 3-2\beta; \xi, \eta) \frac{dx(s)}{ds} ds;$$

$$J_4(x_0, y_0) = \int\limits_{C_\rho} y\left(r^2\right)^{\alpha-\beta-2} F_2(2-\alpha-\beta; 1-\alpha, 1-\beta; 2-2\alpha, 2-2\beta; \xi, \eta) \frac{dy(s)}{ds} ds;$$

$$J_5(x_0, y_0) = \int\limits_{C_\rho} x\left(r^2\right)^{\alpha-\beta-2} F_2(2-\alpha-\beta; 1-\alpha, 1-\beta; 2-2\alpha, 2-2\beta; \xi, \eta) \frac{dx(s)}{ds} ds;$$

$$J_6(x_0, y_0) = \int\limits_0^a x^{2\alpha}\left[y^{2\beta} \frac{\partial q_4(x,y; x_0, y_0)}{\partial n}\right]_{y=0} dx + \int\limits_0^b y^{2\beta}\left[x^{2\alpha} \frac{\partial q_4(x,y; x_0, y_0)}{\partial n}\right]_{x=0} dy.$$

Вводя полярные координаты

$$x = x_0 + \rho\cos\varphi, \quad y = y_0 + \rho\sin\varphi \tag{3.13}$$

в интеграле $J_1(x_0, y_0)$, получим

$$J_1(x_0, y_0) = \int\limits_0^{2\pi} (x_0 + \rho\cos\varphi)(y_0 + \rho\sin\varphi)\left(\rho^2\right)^{\alpha+\beta-2} \times$$
$$\times F_2(3-\alpha-\beta; 1-\alpha, 1-\beta; 2-2\alpha, 2-2\beta; \xi, \eta) d\varphi. \tag{3.14}$$

Исследуем подынтегральное выражение в (3.14). Применяя последовательно известные формулы [11]

$$F_2(a; b_1, b_2; c_1, c_2; x, y) =$$
$$= \sum_{i=0}^{\infty} \frac{(a)_i (b_1)_i (b_2)_i}{(c_1)_i (c_2)_i i!} x^i y^i F(a+i, b_1+i; c_1+i; x) F(a+i, b_2+i; c_2+i; y)$$

и

$$F(a, b; c; x) = (1-x)^{-b} F\left(c-a, b; c; \frac{x}{x-1}\right), \tag{3.15}$$

получим формулу

$$F_2(a; b_1, b_2; c_1, c_2; x, y) = (1-x)^{-b_1}(1-y)^{-b_2} \sum_{i=0}^{\infty} \frac{(a)_i (b_1)_i (b_2)_i}{(c_1)_i (c_2)_i i!}\left(\frac{x}{1-x}\right)^i\left(\frac{y}{1-y}\right)^i \times$$
$$\times F\left(c_1-a, b_1+i; c_1+i; \frac{x}{x-1}\right) F\left(c_2-a, b_2+i; c_2+i; \frac{y}{y-1}\right). \tag{3.16}$$

Воспользовавшись теперь формулой (3.16), функцию Аппеля

$$F_2(3-\alpha-\beta; 1-\alpha, 1-\beta; 2-2\alpha, 2-2\beta; \xi, \eta)$$

запишем в виде

$$F_2(3-\alpha-\beta; 1-\alpha, 1-\beta; 2-2\alpha, 2-2\beta; \xi, \eta) =$$
$$= \left(\rho^2\right)^{2-\alpha-\beta} \left(\rho^2 + 4x_0^2 + 4x_0\rho\cos\varphi\right)^{\alpha-1} \left(\rho^2 + 4y_0^2 + 4y_0\rho\sin\varphi\right)^{\beta-1} P_{11}, \quad (3.17)$$

где

$$P_{11} = \sum_{i=0}^{\infty} \frac{(3-\alpha-\beta)_i (1-\alpha)_i (1-\beta)_i}{(2-2\alpha)_i (2-2\beta)_i \, i!} \times$$
$$\times \left(\frac{4x_0^2 + 4x_0\rho\cos\varphi}{\rho^2 + 4x_0^2 + 4x_0\rho\cos\varphi}\right)^i \left(\frac{4y_0^2 + 4y_0\rho\sin\varphi}{\rho^2 + 4y_0^2 + 4y_0\rho\sin\varphi}\right)^i \times$$
$$\times F\left(-\alpha+\beta-1, 1-\alpha+i; 2-2\alpha+i; \frac{4x_0^2 + 4x_0\rho\cos\varphi}{\rho^2 + 4x_0^2 + 4x_0\rho\cos\varphi}\right) \times$$
$$\times F\left(\alpha-\beta-1, 1-\beta+i; 2-2\beta+i; \frac{4y_0^2 + 4y_0\rho\sin\varphi}{\rho^2 + 4y_0^2 + 4y_0\rho\sin\varphi}\right).$$

Используя известную формулу для $F(a,b;c;1)$ [6]:

$$F(a,b;c;1) = \frac{\Gamma(c)\Gamma(c-a-b)}{\Gamma(c-a)\Gamma(c-b)}, c \neq 0, -1, -2, \ldots, \operatorname{Re}(c-a-b) > 0, \quad (3.18)$$

получим

$$\lim_{\rho \to 0} P_{11} = \frac{\Gamma(2-2\alpha)\Gamma(2-2\beta)}{\Gamma(3-\alpha-\beta)\Gamma(1-\beta)\Gamma(1-\alpha)}. \quad (3.19)$$

Таким образом, согласно (3.14), (3.17) и (3.19), окончательно получим

$$-(2-\alpha-\beta)k_4 x_0^{1-2\alpha} y_0^{1-2\beta} \lim_{\rho \to 0} J_1(x_0, y_0) = -1. \quad (3.20)$$

Далее, учитывая, что

$$\lim_{\rho \to 0} \rho \ln \rho = 0, \quad (3.21)$$

имеем

$$\lim_{\rho \to 0} J_2(x_0, y_0) = \lim_{\rho \to 0} J_3(x_0, y_0) = \lim_{\rho \to 0} J_4(x_0, y_0) = \lim_{\rho \to 0} J_5(x_0, y_0) = 0. \quad (3.22)$$

Наконец, рассмотрим интеграл $J_6(x_0, y_0)$, который, согласно формуле (3.10), можно привести к виду (3.4), т.е.

$$J_6(x_0, y_0) = k(x_0, y_0). \quad (3.23)$$

Теперь, в силу (3.20) – (3.23) из (3.12) следует, что в точке $(x_0, y_0) \in \Omega$ имеет место тождество

$$w_1^{(4)}(x_0, y_0) = k(x_0, y_0) - 1.$$

**Случай 2.** Пусть теперь точка $(x_0, y_0)$ совпадает с некоторой точкой $M_0$, лежащей на кривой Г. Проведем окружность малого радиуса $\rho$ с центром в точке $(x_0, y_0)$. Эта окружность вырежет часть $\Gamma_\rho$ кривой Г. Оставшуюся часть кривой обозначим через $\Gamma - \Gamma_\rho$. Обозначим через $C'_\rho$ часть окружности $C_\rho$, лежащей внутри области $\Omega$, и рассмотрим область $\Omega_\rho$, ограниченную кривыми $\Gamma - \Gamma_\rho$, $C'_\rho$ и отрезками $[0,a]$ и $[0,b]$ осей $x$ и $y$ соответственно. Тогда имеем

$$w_1^{(4)}(x_0, y_0) \equiv \int_0^l x^{2\alpha} y^{2\beta} \frac{\partial q_4(x, y; x_0, y_0)}{\partial n} ds =$$
$$= \lim_{\rho \to 0} \int_{\Gamma - \Gamma_\rho} x^{2\alpha} y^{2\beta} \frac{\partial q_4(x, y; x_0, y_0)}{\partial n} ds. \quad (3.24)$$

Так как точка $(x_0, y_0)$ лежит вне этой области, то в этой области функция $q_4(x, y; x_0, y_0)$ является регулярным решением уравнения (1.1) и в силу (2.4)

$$\int_{\Gamma - \Gamma_\rho} x^{2\alpha} y^{2\beta} \frac{\partial q_4(x, y; x_0, y_0)}{\partial n} ds = \int_0^a x^{2\alpha} \left[ y^{2\beta} \frac{\partial q_4(x, y; x_0, y_0)}{\partial y} \right]_{y=0} dx +$$
$$+ \int_0^b y^{2\beta} \left[ x^{2\alpha} \frac{\partial q_4(x, y; x_0, y_0)}{\partial x} \right]_{x=0} dy + \int_{C'_\rho} x^{2\alpha} y^{2\beta} \frac{\partial}{\partial n} \{q_4(x, y; x_0, y_0)\} ds. \quad (3.25)$$

Подставляя (3.25) в (3.24), с учетом (3.23) и (1.4), получим

$$w_1^{(4)}(x_0, y_0) = k(x_0, y_0) + \lim_{\rho \to 0} \int_{C'_\rho} x^{2\alpha} y^{2\beta} \frac{\partial q_4(x, y; x_0, y_0)}{\partial n} ds. \quad (3.26)$$

Вводя снова полярные координаты (3.13) в интеграле (3.26) и переходя к пределу при $\rho \to 0$, получим

$$\lim_{\rho \to 0} \int_{C'_\rho} x^{2\alpha} y^{2\beta} \frac{\partial q_4(x, y; x_0, y_0)}{\partial n} ds = -\frac{1}{2}.$$

Таким образом,

$$w_1^{(4)}(x_0, y_0) = k(x_0, y_0) - \frac{1}{2}.$$

**Случай 3.** Положим, наконец, что точка $(x_0, y_0)$ лежит вне области $\Omega$. Тогда функция $q_4(x, y; x_0, y_0)$ есть регулярное решение уравнения (1.1) внутри области $\Omega$ с непрерывными производными всех порядков вплоть до контура $\Gamma$ и в силу (2.4)

$$w_1^{(4)}(x_0, y_0) \equiv \int_0^l x^{2\alpha} y^{2\beta} \frac{\partial}{\partial n} \{q_4(x, y; x_0, y_0)\} ds =$$
$$= \int_0^a x^{2\alpha} \left[ y^{2\beta} \frac{\partial q_4(x, y; x_0, y_0)}{\partial y} \right]_{y=0} dx + \int_0^b y^{2\beta} \left[ x^{2\alpha} \frac{\partial q_4(x, y; x_0, y_0)}{\partial x} \right]_{x=0} dy = k(x_0, y_0).$$

Лемма 1 полностью доказана.

**Теорема 1.** *Для любых точек* $(x,y)$ *и* $(x_0, y_0) \in R_+^2$ *при* $x \neq x_0$ *и* $y \neq y_0$, *справедливо неравенство:*

$$|q_4(x,y;x_0,y_0)| \leq C(xx_0)^{1-2\alpha}(yy_0)^{1-2\beta}\left(r_1^2\right)^{\alpha-1}\left(r_2^2\right)^{\beta-1}\ln\left(\frac{r^2}{r_1^2}+\frac{r^2}{r_2^2}-\frac{r^2}{r_1^2}\cdot\frac{r^2}{r_2^2}\right), \quad (3.27)$$

*где* $C$ – *постоянная, а* $\alpha$ *и* $\beta$ – *действительные числа, причем* $0 < 2\alpha, 2\beta < 1$, *а* $r$, $r_1$ *и* $r_2$ – *выражения, определенные в* (1.3).

**Доказательство.** Из (3.16), с учетом неравенства

$$F\left(-\alpha+\beta, 1-\alpha+i; 2-2\alpha+i; 1-\frac{r^2}{r_1^2}\right)F\left(\alpha-\beta, 1-\beta+i; 2-2\beta+i; 1-\frac{r^2}{r_2^2}\right) \leq C_1,$$

получим

$$|q_4(x,y;x_0,y_0)| \leq$$
$$\leq C_1 k_4 \frac{(xx_0)^{1-2\alpha}(yy_0)^{1-2\beta}}{\left(r_1^2\right)^{1-\alpha}\left(r_2^2\right)^{1-\beta}} {}_3F_2\left[\begin{array}{c}2-\alpha-\beta, 1-\alpha, 1-\beta;\\2-2\alpha, 2-2\beta;\end{array}\left(1-\frac{r^2}{r_1^2}\right)\left(1-\frac{r^2}{r_2^2}\right)\right], \quad (3.28)$$

где $C_1 > 0$ – постоянная, а ${}_3F_2\left[\begin{array}{c}a_1,a_2,a_3;\\b_1,b_2;\end{array}z\right]$ – обобщенная гипергеометрическая функция Гаусса [6].

Теперь, согласно формуле [6, 7],

$$F(a_1,a_2,a_3;b_1,b_2;z) = \frac{\Gamma(b_1)\Gamma(b_2)}{\Gamma(a_1)\Gamma(a_2)\Gamma(a_3)}\left\{-\sum_{k=0}^{\infty}\frac{c_k}{k!}(1-z)^k\ln(1-z)+\sum_{k=0}^{\infty}d_k^+(1-z)^k\right\},$$

где $b_1+b_2-a_1-a_2-a_3=0$; $|1-z|<1$; $|\arg(1-z)|<\pi$; $\operatorname{Re} a_j > 0$; $a_j \neq 0, -1, -2, \ldots$;
$j=1,2,3$; $c_k$ и $d_k^+$ известные постоянные, из (3.28) вытекает неравенство (3.27). Теорема 1 доказана.

Таким образом, функция $q_4(x,y;x_0,y_0)$ имеет логарифмическую особенность при $r=0$.

**Теорема 2.** *Если кривая* $\Gamma$ *удовлетворяет перечисленным выше условиям, то*

$$\int_{\Gamma} x^{2\alpha} y^{2\beta}\left|\frac{\partial q_4(x,y;x_0,y_0)}{\partial n}\right|ds \leq C_1,$$

*где* $C_1$ – *постоянная.*

**Доказательство** теоремы 2 следует из условий (3.1) и формулы (3.10).

Формулы (3.3) показывают, что при $\mu_4(s) \equiv 1$ потенциал двойного слоя испытывает разрыв непрерывности, когда точка $(x,y)$ пересекает кривую $\Gamma$. В случае произвольной непрерывной плотности $\mu_4(s)$ имеет место

**Теорема 3.** *Потенциал двойного слоя* $w^{(4)}(x_0,y_0)$ *имеет пределы при стремлении точки* $(x_0,y_0)$ *к точке* $(x(s),y(s))$ *кривой* $\Gamma$ *извне или изнутри. Если предел значений* $w^{(4)}(s)$ *изнутри обозначить через* $w_i^{(4)}(x_0,y_0)$, *а предел извне через* $w_e^{(4)}(x_0,y_0)$, *то для непрерывной плотности* $\mu_4(s)$ *имеют место формулы*

$$w_i^{(4)}(s) = -\frac{1}{2}\mu_4(s) + \int_0^l \mu_4(t) K_4(s,t) dt \qquad (3.29)$$

*и*

$$w_e^{(4)}(s) = \frac{1}{2}\mu_4(s) + \int_0^l \mu_4(t) K_4(s,t) dt, \qquad (3.30)$$

*где*

$$K_4(s,t) = [x(t)]^{2\alpha}[y(t)]^{2\beta} \frac{\partial}{\partial n}\{q_4[x(t),y(t);x_0(s),y_0(s)]\},$$

*точки* $(x(s), y(s))$ *и* $(x_0(t), y_0(t))$ *лежат на кривой* $\Gamma$.

**Доказательство** теоремы 3 следует из леммы 1 и теорем 1 и 2.

Функция

$$w_0^{(4)}(s) = \int_0^l \mu_4(t) K_4(s,t) dt$$

непрерывна при $0 \le s \le l$, что следует из хода доказательства теоремы 3. Принимая во внимание формулы (3.29) – (3.30) и непрерывность функций $w_0^4(s)$ и $\mu_4(s)$ при $0 \le s \le l$, можем утверждать, что потенциал двойного слоя $w^{(4)}(x_0, y_0)$ есть функция, непрерывная внутри области $\Omega$ вплоть до кривой $\Gamma$.

Applying a method of complex analysis (based upon analytic functions), R.P. Gilbert in 1969 constructed an integral representation of solutions of the generalized bi-axially symmetric Helmholtz equation. Fundamental solutions of this equation were constructed recently. In fact, when the spectral parameter is zero, fundamental solutions of the generalized bi-axially symmetric Helmholtz equation can be expressed in terms of Appell's hypergeometric function of two variables of the second kind. All the fundamental solutions of the generalized bi-axially symmetric Helmholtz equation are known, and only for the first one the theory of potential was constructed. In this paper, we aim at constructing a theory of double-layer potentials corresponding to the fourth fundamental solution. Using some properties of Appell's hypergeometric functions of two variables, we prove limiting theorems and derive integral equations containing double-layer potential densities in the kernel.

Keywords: generalized bi-axially symmetric Helmholtz equation; Green's formula; fundamental solution; fourth double-layer potential; Appell's hypergeometric functions of two variables; integral equations with double-layer potential density.

*EHRGASHEV Tuhtasin Gulamzhanovich* (Tashkent Institute of Irrigation and Agricultural Mechanization Engineers, Tashkent, Uzbekistan)
E-mail: ertuhtasin@mail.ru